%% file: ReesAlgebra-JSAG.tex
\documentclass[twoside,12pt, leqno]{amsart}
\usepackage{amsmath,amscd,amsthm,amssymb,amsxtra,latexsym,epsfig,epic,graphics}
\usepackage[matrix,arrow,curve]{xy}
\usepackage{graphicx}
\usepackage{diagrams}
\usepackage{tikz,color}  
\input preamble.tex

\def\PP{{\mathbb P}}

\def\image{{\rm image}}

\def\fix#1{{\bf ***Fix:} #1 {\bf ***}}

\def\RR{{\mathcal R}}

\def\Spec{{{\rm Spec}\,}}

\newarrow{Iso} -----

\def\gr{{\rm gr}}

\def\lbracket{{[\kern-1.5pt[}}
\def\rbracket{{]\kern-1.5pt]}}

\makeatletter
\def\Ddots{\mathinner{\mkern1mu\raise\p@
\vbox{\kern7\p@\hbox{.}}\mkern2mu
\raise4\p@\hbox{.}\mkern2mu\raise7\p@\hbox{.}\mkern1mu}}
\makeatother

\usepackage{times}
\newdimen\x \x=12pt

\usepackage{color}

\usepackage[breaklinks,bookmarksopen,bookmarksnumbered,urlcolor=blue]{hyperref}
\hypersetup{colorlinks=true,backref=true,citecolor=blue}

\author{David Eisenbud}

\title{The Rees Algebra package in Macaulay2}
\begin{document}

\begin{abstract}
This introduces Rees algebras and some of their uses with illustrations via version 2.0 of the Macaulay2 package ReesAlgebra.m2.
\end{abstract}

\maketitle

\section*{Introduction}
A central construction in modern commutative algebra starts from
an ideal $I$ in a commutative ring $R$, and produces the \emph{Rees algebra}
$$
\RR(I) :=  R\oplus I\oplus I^2\oplus I^3\oplus\cdots \cong R[It]\subset R[t]
$$ 
where $R[t]$ denotes the polynomial algebra in one variable $t$ over $R$. For basics on Rees algebras, see \cite{V} and \cite{SW}.

From the point of view of algebraic geometry, the Rees algebra $\RR(I)$ is a homogeneous
coordinate ring for the graph of a rational map whose total space is the blowup
of $\Spec R$ along the scheme defined by $I$.
 (In fact, the  ``Rees Algebra'' is sometimes called the ``blowup algebra''.)

     Rees Algebras were first studied in the algebraic context by
      by David Rees, in the
     now famous paper~\cite{Rees}. Actually
     Rees mainly studied the ring 
     $R[It,t^{-1}]$, now also called the \emph{extended Rees
     Algebra} of I. 
     
Mike Stillman and I, if memory serves, wrote a Rees Algebra script for Macaulay classic. It was augmented, and made into the package ReesAlgebra.m2  around 2002,  to study a generalization of Rees algebras to modules described in \cite{EHU}. Subsequently
Amelia Taylor, 
Sorin Popescu, the present author,
and, at the Macaulay2 Workgroup in July 2017, 
Ilir Dema,
Whitney Liske, and
Zhangchi Chen
contributed routines
for computing many of the invariants of an ideal or module
defined in terms of Rees algebras. This is in fact the package's primary utility, since Rees Algebras
of modules other than ideals are comparatively little studied. 

We first describe the construction and an example from \cite{EHU}. Then we list some of the functionality the package now has; and finally we give an elementary example of how blowups can resolve singularities.

\section{The Rees Algebra of a module} There are several possible ways of extending the Rees algebra construction from ideals to modules. For simplicity we will henceforward only consider finitely generated modules
over Noetherian rings. Huneke and Ulrich and I argued in~\cite{EHU} that the most natural is
to think of $R[It]$ as the image of the map of symmetric algebras
$\Sym(\phi): \Sym_R(I) \to \Sym_R(R) = R[t]$, and to generalize it to the case of
an arbitrary finitely generated module $M$ by setting
$$
\RR(M) = \image \Sym(\phi)
$$ 
where $\phi$ is the \emph{versal} map from $M$ to a free module, defined as 
the composition of 
the diagonal 
embedding 
$$
M \to \oplus_{i=1}^mM.
$$
with the map
$$
\oplus_{i=1}^m\phi_i: \oplus_{i=1}^mM \to R^{m}
$$ 
where
$\phi_1,\dots, \phi_m$ generate $\Hom_R(M,R)$.
It has the property that any map from $M$ to a free module factors (non-uniquely, in general) through $\phi$. 

Though this is not immediate, the Rees algebra of an ideal in any Noetherian ring, in this sense,
 is the same as the Rees algebra in the classical sense, and in most cases one can take
 any embedding of the module into a free module in the definition:
 
\begin{theorem}\label{good cases} [Eisenbud-Huneke-Ulrich, Thms 0.2 and 1.4] Let $R$ be a Noetherian ring
     and let $M$ be a finitely generated $R$-module. Let $\phi: M \to G$
     be the universal embedding of $M$ in a free module, and let
     $\psi: M \to G'$  be any inclusion of $M$ into a free module $G'$. 
     If $R$ is torsion-free over $\ZZ$
     or $R$ is unmixed and generically Gorenstein or $M$ is free locally
     at each associated prime of $R$, or $G=R$, then the image of $\Sym(\phi)$ and the
     image of $\Sym(\psi)$ are naturally isomorphic.
\end{theorem}

Nevertheless some examples do violate the conclusion of Theorem~\ref{good cases}. Here is one from~\cite{EHU}.
     In the following, any finite characteristic would work as well.
     
\begin{footnotesize}
  \begin{verbatim}
p = 5
R = ZZ/p[x,y,z]/(ideal(x^p,y^p)+(ideal(x,y,z))^(p+1))
M = module ideal(z) 
 \end{verbatim}
\end{footnotesize}
\begin{normalsize}
    It is easy to check that $M \cong R^1/(x,y,z)^p$.
     We write $\iota: M\to R^1$ for the embedding as an ideal
     and $\psi$ for the embedding $M \to R^2$ sending $z$ to the vector $(x,y)$.
\end{normalsize}
 \begin{verbatim}
iota = map(R^1,M,matrix{{z}}) 
psi = map(R^2,M,matrix{{x},{y}})
\end{verbatim}
\begin{normalsize}
      Finally, the universal embedding is $M \to R^3$,
     sending $z$ to the vector $(x,y,z)$:
\end{normalsize}
\begin{verbatim}
phi = universalEmbedding(M)
\end{verbatim}
\begin{normalsize}
     We now compute the kernels of the three maps
     on symmetric algebras:
\end{normalsize}
 \begin{verbatim}
Iota = symmetricKernel iota;
Ipsi = symmetricKernel psi;
Iphi = symmetricKernel phi;
\end{verbatim}
\begin{footnotesize}
\begin{normalsize}
      and check that the ones corresponding to $\phi$ and $\iota$
     are equal, whereas the ones corresponding to $\psi$ and $\phi$
     are not---they differ in degree p.
\end{normalsize}
 \begin{verbatim}
i14 : Iiota == Iphi    
o14 = true
i15 : Ipsi == Iphi
o15 = false
i16 : numcols basis(p,Iphi) 
o16 = 3
i17 : numcols basis(p,Ipsi)
o17 = 1
\end{verbatim}
\end{footnotesize}

\section{The Rees Algebra and its relations:\\ {\tt reesIdeal, reesAlgebra, universalEmbedding, symmetricKernel, isLinearType,
associatedGradedRing, normalCone, multiplicity}}

The central routine, reesIdeal (with synonym reesAlgebraIdeal) computes an ideal
defining the Rees algebra $\RR(M)$ as a quotient of a polynomial ring over $R$ from a free presentation of $M$. From the Rees ideal we immediately get
${\tt reesAlgebra}\ M$. In the case when $M$ is an ideal in $R$ we also compute
the important ${\tt associatedGradedRing\ }M = \RR(M)/M$ (and the more geometric sounding but identical ${\tt normalCone\ }M$.)  
If $I$ is a (homogeneous) ideal primary to the
maximal ideal of a standard graded ring $R$ we compute the
Hilbert-Samuel multiplicity of $I$ with ${\tt multiplicity\ }I$.

We now describe the basic computation. Suppose that $M$ has 
a set of generators represented by a map from a free module,
a free presentation
$$
 F\rTo^\alpha M\to 0,
$$
and suppose  $F = R^n$. The symmetric algebra of $F$ over $R$ is then a polynomial ring
$\Sym_R(F_0) = R[t_1,\dots, t_n]$ on $n$ new indeterminates $t_1,\dots, t_n$. By the universal
property of the symmetric algebra there is a canonical surjection
$\Sym_R(F)\rightarrow \Sym_R(M)$, so we may compute the Rees algebra of $M$ as
a quotient of the  $\Sym_R(F)$. The call
$$
I = {\tt reesIdeal}\ M
$$
first calls ${\tt universalEmbedding\ }M$ to compute the versal map from $M$ to a free module $\beta: M\to G$. The call ${\tt symmetricKernel\ } \alpha\circ\beta$  then constructs the map of symmetric algebras $\beta\circ \alpha:\Sym_R(F)\to \Sym_R(G)$ and  calls
 the built-in Macaulay2 routine to compute the kernel 
$$
I = {\tt reesIdeal\ }M = \ker \Sym(\beta\circ\alpha): \Sym_R(F) \to \Sym_R(G).
$$

There is a different way of computing the Rees algebra that is often much more efficient. It begins by constructing the symmetric algebra of $M$, and uses the observation that the construction of the Rees algebra commutes with localization. See \cite[Appendix 2]{E} for the necessary facts about symmetric algebras.

Suppose that $M$ has a free presentation
$$
G\rTo^\alpha F\rTo^\epsilon M\to 0.
$$
The right exactness of the symmetric algebra functor implies that the symmetric algebra of $M$ is the quotient of $\Sym_R(F)$ by an ideal $I_0$ that is
generated by the
entries of the matrix
$$
\begin{pmatrix}
 t_1&\dots&t_n
 \end{pmatrix}
 \circ \phi
$$
(where we have identified $\phi$ with $\Sym_R(F)\otimes_R\phi$).
Thus $I_0$ is generated by polynomials that are linear in the variables $t_i$ (and because
$M$ is the degree 1 part of $\RR(M)$, these are the only linear forms in the $t_i$ in the
Rees ideal.)

If $f\in R$ is an element such that $M[f^{-1}]$ is free on generators $g_1,\dots, g_n$, it follows that after inverting $f$ the Rees algebra of $M$ becomes a polynomial ring over $R[f^{-1}]$ on
indeterminates corresponding to the $g_i$.
$$
\RR(M)[f^{-1}] = \Sym_R(M[f^{-1}]) = R[G_1,\dots, G_n]
$$

 Now suppose in addition that $f$ is a non-zerodivisor in $R$. In the diagram
$$
\begin{diagram}
 \Sym_R(F) &\rTo^\alpha &\Sym_R(M)&\rTo^\beta &\Sym_R(G)\\
 \dTo&&  \dTo&& \dTo\\
 \Sym_R(F)[f^{-1}] &\rTo^\alpha &\Sym_R(M)[f^{-1}]&\rTo^\beta &\Sym_R(G)[f^{-1}]\\
\end{diagram}
 $$
 the two outer vertical maps are inclusions, and it follows that the Rees ideal, which is the
 kernel of the map $\RR(F) = \Sym_R(F) \to \RR(M)$, is equal to the intersection
 of $\RR(F)$ with the kernel of
 $\Sym_R(F)[f^{-1}] \rTo^\beta \Sym_R(G)[f^{-1}]$. This intersection
 may be computed as $I_0:f^\infty$. The call
 $$
 {\tt reesIdeal}(I, f)
 $$
 computes the Rees ideal in this way.
 
More generally, we say that a module $N$ is {\em of linear type} if
the Rees ideal of $M$ is equal to the ideal of the symmetric algebra of $M$; 
for example, any complete intersection ideal is of linear type, and the condition
can be tested by the call
$$
{\tt isLinearType}\ M.
$$
The procedure above really requires only that $f$ be a non-zerodivisor in $R$ and
that $M[f^{-1}]$ be of linear type over $R[f^{-1}]$.

\section{Reductions and the special fiber: {\tt isReduction, minimalReduction, reductionNumber
specialFiber,
specialFiberIdeal,
analyticSpread}}

A \emph{reduction} $J$ of an ideal $I$ is an sub-ideal $J\subset I$ over which $I$ is
\emph{integrally dependent}. In concrete terms this means that there is some integer $r$ such that $JI^r = I^{r+1}$, and the minimal $r$ with this property is called the reduction number.
The property of being a reduction is tested by ${\tt isReduction\ }I$, and the reduction number, is then computed by ${\tt reductionNumber\ } I$. 

Now suppose that $\gm$ is a maximal ideal containing $I$. The special fiber ring is by definition
$\RR(I)/\gm\RR(I)$. It is a standard graded algebra ove the field  $k := R/\gm$, a quotient of 
$\Sym_R(F)/\gm = k[t_1,\dots,t_n]$ where, as before, $F$ is a free module of rank $n$ with a surjection to $M$. The defining ideal of the special fiber ring, and the ring itself, are computed  by the calls
${\tt specialFiberIdeal\ }I$ and ${\tt specialFiberRing\ }I$. 

The dimension of the special fiber ring is called the analytic spread of $I$, usually
denoted
$$
\ell(I) = {\tt analyticSpread\ }I.
$$
Northcott and Rees~\cite{NR} proved that if $k$ is infinite then there always exist reductions
generated by $\ell(I)$ elements, and this is the minimum possible number; these are called
minimal reductions. The smallest possible reduction number for $I$ with respect to an minimal reduction is  by definition ${\tt reduction number} I$ (this is always achieved by any ideal generated by $\ell(I)$ sufficiently general scalar linear combinations of the generators of $I$; but note that when $I$ is homogeneous but has generators of different degrees such linear combinations are sometimes necessarily inhomogeneous.)

An interesting special case occurs when $R$ is a graded ring over $k = R_0$ and the generators $g_1,\dots, g_n$ of $I$ are all homogeneous of the same degree. In this case the special fiber ring is easily seen to be equal to the subring $k[g_1,\dots,g_n]$ (usually \emph{not} a polynomial ring) generated by the elements $g_i$.

\def\G{{\mathbb G}}
For example, if $I$ is the ideal of $p\times p$ minors of a $p\times (p+q)$ matrix, then
the special fiber ring is equal to the homogeneous coordinate ring $\G$ of the Grassmannian of
$p$-planes in $p+q$ space. It follows that $\ell(I) = \dim \G = pq+1$, and the reduction number of $I$ is
$(p-1)(q-1)$.

\section{Finding elements of the Rees ideal: {\tt jacobianDual, expectedReesIdeal}}

Let $M$ be an $R$-module and let $\phi: R^{s}\to R^{m}$ be its presentation matrix.
     The symmetric algebra of I has the form      
     $$
     \Sym_R(I) = Sym(R^{m)}/(T\phi)
     $$
     where we have written $(T\phi)$  for the ideal generated by the entries
     of the product 
 $$
\begin{pmatrix}
 T_{0}&\dots&T_{m}
\end{pmatrix}\phi
$$
and the $T_{i}$ correspond to the generators of $I$. If 
     $$
     X = \begin{pmatrix}
x_1&\dots&x_{n}
\end{pmatrix}
$$
     with $x_i \in R$, and the ideal $J$ generated by the entries of $X$ 
     contains the entries of the matrix $\phi$, then there is 
     a matrix $\psi$ defined over $R[T_0..T_m]$, called the Jacobian Dual of $\phi$ with respect to $X$,
     such that $T\phi = X\psi$. (the matrix $\psi$ is generally
     not unique; Macaulay2 computes it using Gr\"obner division with remainder.)
           
     If $I,J$ each contain a non-zerodivisor then
     $J$ will have grade $\geq 1$ on the Rees algebra $\RR(I)$. Since $(T\phi)$ is contained in the
     defining ideal of the Rees algebra, the vector $X$ is annihilated by the matrix
     $\psi$ when regarded over the Rees algebra, and the relation
     $X\psi \equiv 0$ in $\RR(I)$ implies that the $m\times m$ minors of $\psi$ are
     in the Rees ideal of $I$. In very favorable circumstances,
     one may even have the equality 
     {\tt reesIdeal I = ideal(T*phi)+ideal minors(psi)}. We illustrate with a Theorem of
     Morey and Ulrich. Recall that an ideal $I$ is said to satisfy the condition
     $G_{m}$ if the number of generators of the localized ideal $I_{P}$ is $\leq \codim P$
     for every prime ideal $P$ of codimension $<m$; equivalently, if $I$ has  presentation
      matrix $\phi$ as above, 
     $$
     \codim I_{m-p}(\phi)>p
     $$
     for $1<= p < \ell$

\begin{theorem}[\cite{MU}]
    Let $R$ be a local Gorenstein ring with infinite residue field, let $I$ be a perfect ideal
     of grade 2 with m generators, and let $\phi$ be the presentation matrix of $I$,
     and let $\psi$ be the Jacobian dual matrix. Let
     $\ell = \ell(I)$ be the analytic spread. Suppose that
     $I$ satisfies the condition $G_{\ell}$. The following conditions are equivalent:
 
\begin{enumerate}
     \item $\RR(I)$ is Cohen-Macaulay and $I_(m-\ell)(phi) = I_1(phi)^{m-\ell}$.
     \item $r(I) < \ell$ and $I_{m+1-\ell}\phi = (I_1\phi)^{m+1-\ell}$.
     \item The ideal of $\RR(I)$ is equal to the sum of the ideal
     of $\Sym(I)$ with the Jacobian dual minors, $I_{m}\psi$.)
\end{enumerate}
    
\end{theorem}

We can check all these conditions with functions in the package.      
     We start with the presentation matrix phi of an $m = n+1$-generator perfect ideal
     Such that the first row consists of the $n$
     variables of the ring, and the rest of whose rows are reasonably general (in this
     case random quadrics):

\begin{footnotesize}
 \begin{verbatim}
i1 : loadPackage("ReesAlgebra", Reload =>true)
i2 : setRandomSeed 0
i3 : n=3;
i4 : kk = ZZ/101;
i5 : S = kk[a_0..a_(n-2)];
i6 : phi' = map(S^(n),S^(n-1), (i,j) -> if i == 0 then a_j else random(2,S));
             3       2
o6 : Matrix S  <--- S
i7 : I = minors(n-1,phi');
\end{verbatim}
\end{footnotesize}
This is a perfect, codimension 2 ideal, as we see from the Betti table:
\begin{footnotesize}
 \begin{verbatim}
i8 : betti (F = res I)
            0 1 2
o8 = total: 1 3 2
         0: 1 . .
         1: . . .
         2: . 2 .
         3: . 1 2
\end{verbatim}
\end{footnotesize}
As we constructed the matrix {\tt phi'} it was not homogeneous, but the resolution is, so we take instead:
\begin{footnotesize}
 \begin{verbatim}
i9 : phi = F.dd_2;
 \end{verbatim}
 \end{footnotesize}
We  compute the analytic spread $\ell$ and the reduction number $r$:
\begin{footnotesize}
 \begin{verbatim}
i12 : ell = analyticSpread I
i13 : r = reductionNumber(I, minimalReduction I)
o13 = 1
\end{verbatim}
\end{footnotesize}
Now we can check the condition $G_{ell}$, first probabilistically:
\begin{footnotesize}
 \begin{verbatim}
i15 : whichGm I >= ell
o15 = true
\end{verbatim}
\end{footnotesize}
and now deterministically:
\begin{footnotesize}
 \begin{verbatim}
i17 : apply(toList(1..ell-1), p-> {p+1, codim minors(n-p, phi)})
o17 = {{2, 2}}
\end{verbatim}
\end{footnotesize}
We now check the three equivalent conditions of the Morey-Ulrich Theorem.
           Note that since $\ell = n-1$ in this case, the second part of conditions
           1) and 2) is vacuously satisfied, and since $r<\ell$
           the conditions must all be satisfied.
           We first check that $\RR(I)$ is Cohen-Macaulay:
\begin{footnotesize}
 \begin{verbatim}
i19 : reesI = reesIdeal I;
o19 : Ideal of S[w , w , w ]
                  0   1   2
i20 : codim reesI
o20 = 2
i21 : betti res reesI
             0 1 2
o21 = total: 1 3 2
          0: 1 . .
          1: . . .
          2: . 2 .
          3: . 1 2
\end{verbatim}
\end{footnotesize}
Finally, we wish to see that reesIdeal I is generated by the ideal 
           of the symmetric algebra together with the jacobian dual:
\begin{footnotesize}
 \begin{verbatim}
i23 : psi = jacobianDual phi
o23 = {0, 1} | 11w_1a_0-37w_1a_1                 w_2                              |
      {0, 1} | -11w_0a_0-34w_0a_1-46w_1a_1+14w_2 11w_0a_1-11w_1a_0-30w_1a_1+18w_2 |
                            2                     2
o23 : Matrix (S[w , w , w ])  <--- (S[w , w , w ])
                 0   1   2             0   1   2
\end{verbatim}
\end{footnotesize}
We now compute the ideal J of the symmetric algebra; the call {\tt symmetricAlgebra I}
           would return the ideal over a different ring, so we do it by hand:
\begin{footnotesize}
 \begin{verbatim}
i25 : ST = ring psi
i26 : T = vars ST
o26 = | w_0 w_1 w_2 |
i27 : J = ideal(T*promote(phi, ST))
                            2          2               2                 2                    2
o27 = ideal ((- 11a a  - 34a )w  + (11a  - 37a a  - 46a )w  + 14a w , 11a w  + (- 11a a  - 30a )w  + (a  +
                   0 1      1  0       0      0 1      1  1      1 2     1 0         0 1      1  1     0  
      --------------------------------------------------------------------------------------------------------
      18a )w )
         1  2
i28 :      betti res J
             0 1 2
o28 = total: 1 2 1
          0: 1 . .
          1: . . .
          2: . 2 .
          3: . . .
          4: . . 1
i29 : J1 = minors(ell, psi)
                        2              2             2  2                                               2
o29 = ideal((20a a  - 3a )w w  + (- 20a  - 24a a  - a )w  + (11a  + 34a )w w  + (- 4a  - 14a )w w  - 14w )
                0 1     1  0 1         0      0 1    1  1       0      1  0 2        0      1  1 2      2
  \end{verbatim}
  \end{footnotesize}
And we compute the resolution of $J+J1$, to see that the resulting ideal is perfect, which also shows that
it is the full ideal of the Rees algebra. We also check directly that it has the same resolution as the computed
Rees ideal of $I$:
\begin{footnotesize}
 \begin{verbatim}
i30 : betti (G = res trim (J+J1))
             0 1 2
o30 = total: 1 3 2
          0: 1 . .
          1: . . .
          2: . 2 .
          3: . 1 2
i31 : betti res reesIdeal I
             0 1 2
o31 = total: 1 3 2
          0: 1 . .
          1: . . .
          2: . 2 .
          3: . 1 2
o31 : BettiTally
  \end{verbatim}
  \end{footnotesize}

\section{Distinguished subvarieties: {\tt distinguished, intersectInP}}

The key construction in the Fulton-MacPherson definition of the refined intersection  product 
\cite[Section 6.1]{F}
involves normal cones, and is easy to implement using the tools in this package. The simplest case is the intersection of two subvarieties $X,V$. If $X$ and $V$ meet in the \emph{expected dimension}, defined to be $\dim V - \codim_{Y} X$,
and the ambient variety $Y$ is smooth, then one can assign multiplicities to the components $W_{i}$ of $X\cap V$, and the $X\cdot V$ is a positive linear
combination of these components. The astonishing result of the theory is that if $X\subset Y$ is locally a complete intersection, then, no matter how singular $Y$ and no matter how strange the actual intersection $X\cap V$, the intersection product 
$X\cdot V$ can be given a meaning as rational equivalence class of cycles of the expected dimension on $X$, or even on $X\cap V$. This class comes with a canonical decomposition $\sum_{i}m_{i}\alpha_{i}$, 
where the $m_{i}$ are positive integers, and $\alpha_{i}$ is a cycle of the expected dimension (possibly 0)
on a certain subvariety $Z_{i}\subset X\cap V$ called a distinguished variety of the intersection (the $Z_{i}$ need not be distinct.)

In the general case, the subvariety $V$ is replaced by a morphism $f:V\to Y$ from a variety $V$, and this is the key to the functoriality of the intersection product. The routines in this package work in the general setting, but for simplicity we will stick with the basic case in this description.

We now describe the distinguished subvarieties and their multiplicities. This part of the construction sheafifies, so (as in the package) we work in the affine case. We do not require any hypothesis on $X, Y$ or $V$. 

Let $S$ be the coordinate ring of $Y$ and let $I\subset S$ be  ideal of $X$
Write
$
\gr_{I}S
$ for the associated graded ring $S/I\oplus I/I^{2}\oplus\dots$ of $I$ in $S$, and let
$\pi$ be the inclusion of $S/I$ into $T$ as the degree 0 part.

Let $R$ be the coordinate ring of $V$, and let $f:S\to R$ be a morphism (if $V$ is a subvariety of $R$ then $f$ will be a projection  $S\to S/J$. Let $K\subset T$ be the kernel of the induced map $ \gr_{I}S \to \gr_{f(I)}R$.  

Let $P_{1},\dots, P_{m}$ be the minimal primes over $K$ in $\gr_{I}R$We define $p_{i}$ to be the degree 0 part of $P_{i}$, that is, $p_{i} := P_{i}\cap S/I$. These are the distinguished prime ideals of $S/I$, and they clearly contain the kernel of $\overline f: S/I \to R/f(I)$, so in the case where $R = S/J$ they contain $I+J$, so they really are subvarieties
of $X\cap V$.

We further define the multiplicity $m_{i}$ to be the multiplicity with which $P_{i}$ appears in the primary decomposition of $K$ -- that is,
$$
m_{i} := \length_{T_{P_{i}}}{P_{i}}_{P_{i}}/K_{P_{i}}.
$$

Returning to the geometric language, and the case where $X\subset Y$ is locally a complete intersection in a quasi-projective variety, the cycle class $\alpha_{i}$ in the Chow group of the variety 
$Z_{i}$ corresponding to $p_{i}$ is defined as the Gysin image of the class of the subvariety corresponding
to $P_{i}$ in the projectivized normal bundle of $X$ in $Y$---a construction not included in this package. 

Here are some simple examples in which {\tt distinguished} is used to compute the distinguished varieties of
intersections in ${\mathbb A}^{n}$, via the function {\tt intersectInP}. First, the intersection of a conic with a tangent line.
\begin{footnotesize}
 \begin{verbatim}
i2 : kk = ZZ/101;
i3 : P = kk[x,y];
i4 : I = ideal"x2-y";J=ideal y;
i6 : intersectInP(I,J)
o6 = {{2, ideal (y, x)}}
\end{verbatim}
\end{footnotesize}
Slightly more interesting, the following shows what happens when the intersections aren't rational:
\begin{footnotesize}
 \begin{verbatim}
i7 : I = ideal"x4+y3+1";
i8 : intersectInP(I,J)
                     2                        2
o8 = {{1, ideal (y, x  + 10)}, {1, ideal (y, x  - 10)}}
\end{verbatim}
\end{footnotesize}
The real interest in the construction is in the case of improper intersections. Here are some typical results:{footnotesize}
\begin{footnotesize}
 \begin{verbatim}
i9 : I = ideal"x2y";J=ideal"xy2";
i11 : intersectInP(I,J)
o11 = {{2, ideal x}, {5, ideal (y, x)}, {2, ideal y}}
i12 : intersectInP(I,I)
o12 = {{1, ideal y}, {4, ideal x}, {4, ideal (y, x)}}
\end{verbatim}
\end{footnotesize}

\section{Rees Algebras and Desingularization}

We conclude this note with an example illustrating a general result about projective birational maps of varieites.
Recall that a map $B\to X$ of varieties is projective if it is the composition of a closed embedding $B\subset X\times \PP^n$ with the projection to $X$.
It is birational if it is generically an isomorphism. The inclusion of a ring into the Rees algebra
of an ideal corresponds to a map from Proj of the Rees algebra to spec of the ring, called a blowup, that is such a proper birational transformation, and in fact every proper birational transformation to an affine variety (or more generally to any scheme, if one works with sheaves of ideals) can be realized in this way.

The Theorem of embedded resolution of singularities, proven by Hironaka in characteristic 0 and conjectured in general, says that give any subvariety $X$ of a smooth variety $Y$, there is
a finite sequence of blowups 
$$
B_n \to \cdots B_2 \to B_1 \to Y
$$
of smooth subvarieties  and a component of the preimage of
$X$ in $B_n$ that is smooth. In the case of plane curves, this can be done with a sequence of blowups of closed points. But in fact *any* sequence of blowups of a quasi-projective variety can be replaced with a single blowup (\cite[Theorem II.7.17]{Hartshorne} of a more complicated ideal. We illustrate with the desingularization of a tacnode (the union of two smooth curves that meet with a simple tangency.)

\begin{example}
Blowing-up $(x^2,y)$ in k[x,y] desingularlizes the tacnode $x^2-y^4$ in a single step. 
\end{example}
\begin{footnotesize}
 \begin{verbatim}
 Macaulay2, version 1.10
with packages: ConwayPolynomials, Elimination, IntegralClosure, InverseSystems, LLLBases, PrimaryDecomposition, ReesAlgebra, TangentCone
i1 : R = ZZ/32003[x,y];
i2 : tacnode = ideal(x^2-y^4);
i3 : mm = ideal(x,y^2);
i4 : B = first flattenRing reesAlgebra mm;
i5 : irrelB = ideal(w_0,w_1);
i6 : proj = map(B,R,{x,y});
i7 : totalTransform = proj tacnode
              4    2
o7 = ideal(- y  + x )
i8 : netList (D = decompose totalTransform)
     +-----------------------+
o8 = |ideal (y, x)           |
     +-----------------------+
     |        2              |
     |ideal (y  + x, w  + w )|
     |                0    1 |
     +-----------------------+
     |        2              |
     |ideal (y  - x, w  - w )|
     |                0    1 |
     +-----------------------+
i9 : exceptional = proj mm
                2
o9 = ideal (x, y )
i10 : strictTransform = saturate(totalTransform, exceptional);
i11 : netList decompose strictTransform
      +-----------------------+
      |        2              |
o11 = |ideal (y  + x, w  + w )|
      |                0    1 |
      +-----------------------+
      |        2              |
      |ideal (y  - x, w  - w )|
      |                0    1 |
      +-----------------------+
i12 : sing0 = sub(ideal singularLocus strictTransform, B);
i13 : sing = saturate(sing0,irrelB)
o13 = ideal 1
\end{verbatim}
\end{footnotesize}
The last line asserts that the singular locus of the the variety ``properTransform'' is empty;
that is, the the scheme defined by "properTransform" is smooth (in this case it is the union
of two disjoint smooth curves.)

\let\thefootnote\relax\footnote{
\noindent A\S Subject Classification:\\
Primary: 13C40, 13H10, 14M06, 14M10;
Secondary: 13D02 , 13N05, 14B12, 14M12.\smallbreak
The author is grateful to the
National Science Foundation for partial support.<}
\bibliographystyle{alpha}

\bigskip

\vbox{\noindent Author Addresses:\par
\smallskip
\noindent{David Eisenbud}\par
\noindent{Mathematical Sciences Research Institute,
Berkeley, CA 94720, USA}\par
\noindent{de@msri.org}\par
}

\end{document}

%% file: preamble.tex
\def\antiddot{\mathinner{\mkern1mu\raise1pt\vbox{\kern7pt\hbox{.}}\mkern2mu
        \raise4pt\hbox{.}\mkern2mu\raise7pt\hbox{.}\mkern1mu}}

\newcommand{\PP}{{\mathbb P}}

\newcommand{\RR}{{\mathbb R}}

\newcommand{\ZZ}{{\mathbb Z}}

\newcommand{\punkt}{\hspace{-.3ex}\raise.15ex\hbox to1ex{\Huge.}}

\def \fix#1 {{\hfill\break \bf (( #1 ))\hfill\break}}

\DeclareMathOperator{\Sym}{Sym}

\DeclareMathOperator{\Spec}{Spec}

\DeclareMathOperator{\Hom}{Hom}

\DeclareMathOperator{\length}{length}

\DeclareMathOperator{\codim}{codim}

\newcommand{\gm}{\mathfrak m}

\def\gr{{\mathfrak {gr}}}

\newtheorem{theorem}{Theorem}[section]

\theoremstyle{definition}

\newtheorem{example}[theorem]{Example}